\documentclass[12pt, a4paper]{amsart}

\usepackage{amsmath}
\usepackage{amssymb}
\usepackage{amsthm}
\usepackage{graphicx}

\begin{document}

\author{Ron Larham*}
\address{*BAE Systems, Underwater Systems Div, Elletra Av, Waterlooville, Hants, PO7 7XS, UK}
\email{ronald.larham@baesystems.com}

\title{Validation of a Model of the Domino Effect?}

\begin{abstract}
A recent paper proposing a model of the limiting speed of the domino effect is discussed with reference to its need and the need of models in general for validation against experimental data. It is shown that the proposed model diverges significantly from experimentally derived speed estimates over a significant range of domino spacing using data from the existing literature and this author's own measurements, hence if its use had had economic importance its use outside its range of validity could have led to loses of one sort or another to its users.
\\
\\
\end{abstract}

\maketitle


\pagestyle{myheadings}
\thispagestyle{plain}
\markboth{R. Larham}{Domino Model Validation}

\section{Introduction}
The purpose of mathematical modelling is to capture 
some aspect of the behaviour of reality in a set of
equations, or some other mathematical structure, to
allow us to examine its behaviour in the comfort of
our offices. But any mathematical model introduces some level
of simplification so before we place our confidence
in a model it is a good idea to compare the behaviour that the
model predicts with reality at least for some subset
of cases of interest. This is validation, and before we
use the predictions of a model we should have shown
that at the very least it reproduce the known behaviour
of the modelled system. Any student of mathematical modelling
must learn the importance of validation if their work is
to be applied with any confidence.

In this note I investigate the validity of a model of
domino wave speed propose by Efthimiou and Johnson \cite{Efth}.
The reason that I am focusing on this model is that
in their paper the authors present their model and
derive the speed in terms of a domino spacing parameter
($d/H$ the ratio of the gap between adjacent dominoes
and their height, see figure \ref{geom}), but do not quote any
experimental results in to support the model despite such data being
easily acquired.

\section{Background}
A number of people have previously constructed models to represent
the propagation of a wave of domino collapse. Some of these
are the models of van Leeuwen \cite{vanLee}, McLachlan et al.
\cite{McLach} and   Banks \cite{Banks} which are references
in \cite{Efth}, and Stronge \cite{Stronge} and Stronge
and Shu \cite{StrongeShu}.
Of these McLachlan et al's model is for "thin" dominoes (or for geometrical
similar dominoes that is the height to thickness ratio is constant),
and Bank's models assume that only pairs of adjacent dominoes are involved
at any one time in the collapse wave. Efthimiou and Johnson's model makes
both of these assumptions.
The others model "thick" dominoes (that is of non-negligible thickness) and
also van Leeuwen and Stronge and Shu model the effects of multiple dominoes
being involved in the wave
(it seems that this is also implicit in McLachlan et al)) they also
model the interactions between the dominoes in more detail.
(McLachlan et al's use a dimensional argument which is independent of any
assumption about the number of dominoes involved in the interactions
giving rise to the collapse wave)
All of these quote some experimental results
though \cite{vanLee} and \cite{Banks} quote the experimental
results of \cite{McLach}. So there is no shortage of data
against which a model may be validated.

My first reaction when I read \cite{Efth}, before I had chased up
references, was that it must be easy to measure the wave speed.
As a result I set about collecting equipment to make my own
measurements, which I will compare with \cite{Efth}'s theoretical
calculations and other experimental data below.

\section{Theory}
McLachlan et al. \cite{McLach} conclude from dimensional analysis that the limiting
wave speed $V$  for thin dominoes satisfies:
\begin{eqnarray*}
 \frac{v}{\sqrt{gH}} = G(L/H) 
\end{eqnarray*}
for some function $G$. Which for thin dominoes is the same as:
\begin{eqnarray*}
 \frac{v}{\sqrt{gH}} = G_1(d/H) 
\end{eqnarray*}

\begin{figure}[ht]
\begin{center}
\includegraphics[width=5.0 in]{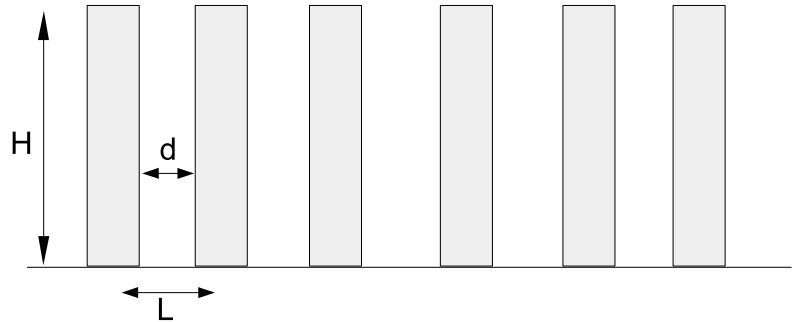}
\end{center}
\caption{Geometry of Domino Array} 
\label{geom} 
\end{figure}

Where $H$ is the height of the dominoes, $d$ the gap between adjacent
dominoes, and $L$ the distance between equivalent points on neighbouring
dominoes (that is the pitch of the domino array) The proposed model of
Efthimiou and Johnson aims to determine what the form of the function
$G_1$ must be (see figure \ref{geom} for the significance of the variables).
\\

Notes: There are additional invisible arguments in functions $G$ and $G_1$ as
the normalised speed may also depend on the dimensionless constants
relevant to the system, in this case these include the coefficient of
friction between dominoes, and the coefficient of restitution for inter
domino impacts. The coefficient of friction between the surface and the dominoes
is of lesser relevance as in domino experiments it is usual to arrange things so that
there is no slipping between the dominoes and the surface. Efthimiou and Johnson
make particular assumptions which mean that none of these dimensionless parameter
appear explicitly in their model.

\section{Generation of Validation Data}
Initially I had toyed with the idea of videoing domino waves, then
extracting the speed from an analysis of the video's frames. I 
abandoned this approach when I realised that audio recording would be
more convenient. The way that I decided to measure the domino wave speed
was to use the sound recorder and microphone on an old laptop to
record the sound of a domino array collapse. (This is far less demanding
in terms of cost of equipment than the high speed photography reported
in \cite{Stronge} and \cite{StrongeShu}). Then to analyse the recording
to extract the frequency of dominoes hitting the their adjacent domino.
This signal is encoded in the envelope of the recording so analysis
techniques analogous to the processing in a crystal AM radio
receiver, or a simple form of DEMON (Detection of Envelope Modulation
On Noise) analysis  similar to that used in passive Sonar processing is required.
The initial sections of each recording were progressively discarded to
identify and eliminate any start up transients. For most of the recordings
the transients were at most slight and easily eliminated but four (the
two with the closest and the two with the widest relatively spacing of
the  dominoes) must be regarded with caution as the results were irregular
for these (they should probably be repeated more carefully with a better
surface and more regular dominoes than the high street ones that I used).

The analysis of the acoustic recordings is fascinating problem in
itself but I will not  go into the detail of the experiment,
it is the results and their relationship to the model predictions 
that I am interested in here. A detailed description of the experimental
method can be found in \cite{RonL1} and \cite{RonL2}.

\section{Results}
\begin{figure}[ht]
\begin{center}
\includegraphics[width=5.0 in]{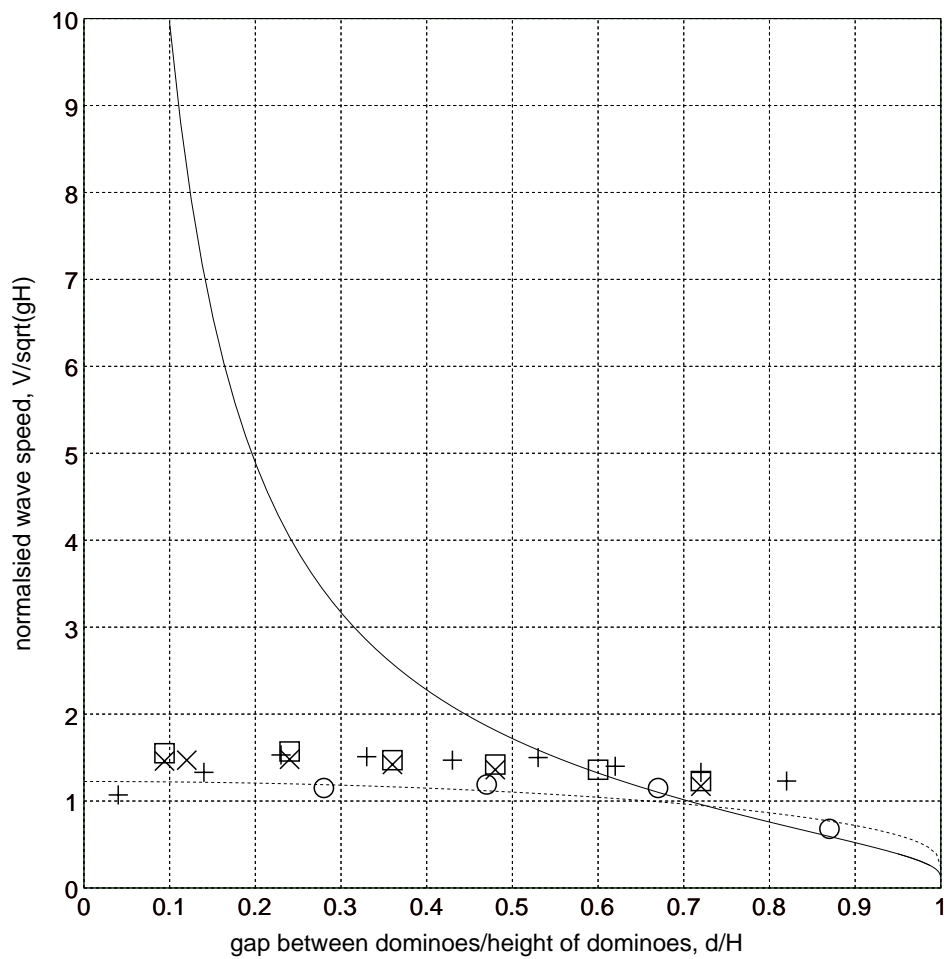}
\end{center}
\caption{Domino Wave Speed Plot; the solid line is the model prediction from
Efthimiou and Johnson \cite{Efth}, dashed line from Banks \cite{Banks},
$+$ measurements by the current author with the dominoes vertical, $\circ$
with them horizontal, $\times$ and $\Box$ measurements from Strong and Shu
\cite{StrongeShu}} 
\label{graph1} 
\end{figure}

All of the papers that report domino wave speed measurements
report speeds $\sim 0.5$ to $1.7 \times \sqrt{gH}$. These are in
broad agreement with my own measurements shown in figure \ref{graph1}.

Given the usual shapes of dominoes I would hope that the thin domino
approximation would be not un-reasonable down to values of $d/H \sim 0.3$. 

As can also be seen in the figure the measured data are comparable to
the predictions of \cite{Efth} over a rather limited range of
$d/H$. This is in contradistinction to the models in the references
(the predictions of Bank's \cite{Banks} is reproduced in the figure for
reference) which in general give rather better agreement with experiment.
The models which represent the effects of multiple dominoes being involved in the
collapse wave being rather better than Bank's model. Even so the
reasonable agreement between the experimental data and the model
predictions from Banks \cite{Banks} is worth noting as it indicates that
the single neighbour domino interaction assumption is not entirely misleading.

\section{Discussion}
The experimental data may be summarised as telling us that to a fair
(hand-waving) approximation for common dominoes the normalised wave speed is independent
of the normalised inter domino interval for all practical intervals. Also that the normalised wave speed is of the order 1, with variation
about this value probably explained by the variation of sliding friction
and possibly coefficient of restitution between different species of dominoes
(see \cite{StrongeShu} and \cite{vanLee}).

\begin{figure}[ht]
\begin{center}
\includegraphics[width=5.0 in]{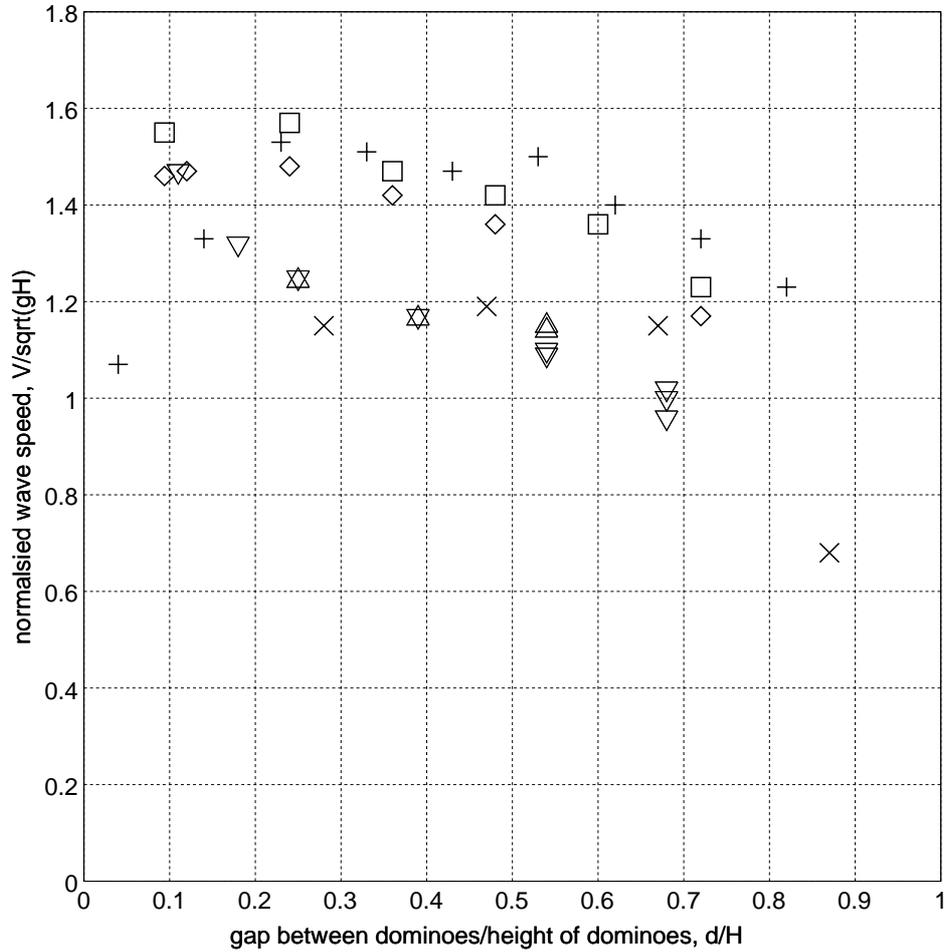}
\end{center}
\caption{Domino Wave Speed Data;
$+$ measurements by the current author with the dominoes vertical, $\times$
with them horizontal, $\square$ and $\lozenge$ measurements from Strong and Shu 
\cite{StrongeShu}, $\triangledown$ (normal dominoes) and $\vartriangle$ (double
height dominoes) measurements from McLachlan et al.
\cite{McLach} } 
\label{graph2} 
\end{figure}

It is difficult to see the variability in the experientially determined
wave speeds in \ref{graph1} due to the vertical scale required to show
the model predictions. The experimental data shown in figure \ref{graph1}
together with that from \cite{McLach} are shown in figure \ref{graph2}
with a vertical scale more suited to showing this data. In this plot the
variation of wave speed with spacing parameter can be seen. For most of the data sets there appears to be a weak dependence of normalised wave
speed on  the spacing parameter $h/H$. The data form McLachlan et al
\cite{McLach} shows what appears to be a stronger dependence on the spacing parameter apparently rising as the spacing parameter becomes small.
Given the other data and McLachlan et al's description of their experimental
method (hand timing the collapse of an array of $\sim ~100$ dominoes)
this is possibly a misleading trend. 

From figure \ref{graph2} we can see that all the reliable data points
give normalised wave speeds in the range $\sim 1 - 1.6$.

(Note the apparent systematic upward curve in McLachlan et al's data
as the domino spacing falls may be due to systematic effects in the
experimental method that was employed.)

\section{Summary}
From the comparison of the model of Efthimiou and Johnson \cite{Efth} and experiment
we see that the area of agreement of experiment and model is rather limited.
Had the model been part of a project with some economic impact we would have been
at risk of being found to not have shown due diligence, which could result in
unfavourable consequences for us and/or our employers in the event of
a failure.

Validation of models is not a chore that we may do after the interesting
parts of a study are completed but an essential activity if our work
is not to be nugatory.

It is also worth while comparing the predictions in the literature with ones
current models predictions, the differences may be important and in need of
explanation

\section*{Acknowledgements}
The author thanks his partner and children for putting up
with domino experiments on the kitchen table extending over
several evenings without expressing any more than slight
derision (or interest), and cats for not prematurely disturbing 
the domino arrays, and Bill Stronge for his advice on some aspects
of the content of this note. Also I want to acknowledge my debt the giants in
who's footprints I stand.

\begin{verbatim}
  
  
\end{verbatim}

\appendix
\begin{center}
\section{Experimental Data}
\end{center}

\begin{verbatim}
  
\end{verbatim}

For reference I include here the data from my experiments on domino wave speed.
\\

\begin{table}[!h]
\begin{center}
\caption {Experimental Results With Dominoes Vertical (italic script indicates less reliable data)}
\begin{tabular}{|c|c c c c c c c c c|}
\hline
$d/H$ & 0.04 & 0.14 & 0.23 & 0.33 & 0.43 & 0.53 & 0.62 & 0.72 & 0.82\\
\hline
$V/\sqrt{gH}$ & \sl{1.07} & \sl{1.33} & 1.53 & 1.51 & 1.47 & 1.50 & 1.40 & 1.33 & \sl{1.23}\\
\hline
\end{tabular}
\label{table1}
\end{center}
\end{table}

\begin{table}[!h]
\begin{center}
\caption {Experimental Results With Dominoes Horizontal (italic script indicates less reliable data)}
\begin{tabular}{|c|c c c c|}
\hline
$d/H$ & 0.28 & 0.47 & 0.67 & 0.87\\
\hline
$V/\sqrt{gH}$ & 1.15 & 1.19 & 1.15 & \sl{0.68}\\
\hline
\end{tabular}
\label{table2}
\end{center}
\end{table}

\textbf{Notes}
The last entry in table \ref{table2} has a spacing greater than the maximum
for which one would expect the domino wave to propagate. At a value of $d/H > \sqrt{3}/2$
a domino strikes its neighbour below its' mid point, and under these conditions
it may well not topple in the expected manner, this is van Leeuwen's practical upper limit
for the wave to propagate. So it is no surprise that
the data for this point is unreliable and this was the largest spacing at which 
I could get the wave to propagate. Presumably it did propagate in this case
as a result of the irregularities in the domino geometry and spacing, or some
other unidentified reason.

\end{document}